\documentclass[reqno]{amsart}
\usepackage{amssymb,latexsym}
\usepackage{txfonts}
\usepackage{amscd}
\usepackage{amsmath}
\usepackage{graphicx}

\def\ba{\begin{array}}
\def\ea{\end{array}}
\def\be{\begin{equation}}
\def\ee{\end{equation}}
\def\bee{\begin{eqnarray}}
\def\beee{\begin{eqnarray*}}
\def\eee{\end{eqnarray}}
\def\eeee{\end{eqnarray*}}
\def\nn{\nonumber}

\numberwithin{equation}{section}

\newcommand{\e}{\overline{\varepsilon}}

\newtheorem{thm}{Theorem}[section]
\newtheorem{lem}{Lemma}[section]

\newtheorem{rem}{Remark}[section]
\newtheorem{pro}{Proposition}[section]
\newtheorem{defi}{Definition}[section]

\newcommand{\noi}{\noindent}

\def \pr {\noindent {\it Proof.} }

\def \eop {  \hfill $\Box$}

\def\R{\Bbb R}

\def \a {\alpha}
\def \b {\beta}

\def \e {\epsilon}

\def \h {\hspace{.5cm}}

\def \ra {\longrightarrow}

\def \h* {\hspace*{1cm}}
\def\la{\langle}
\def\ra{\rangle}

\title[Regularity for harmonic maps into certain Pseudo-Riemannian manifolds]{Regularity for harmonic maps into certain Pseudo-Riemannian manifolds}

\author{Miaomiao Zhu}

\address{Mathematics Institute, University of Warwick \\  CV4 7AL, Coventry, UK}
\email{Miaomiao.Zhu@warwick.ac.uk}
\subjclass[2000]{58E20; 53C50; 35J60; 35B65}
\keywords{ Harmonic map; Regularity; Lorentzian manifold; Pseudo-Riemannian manifold}
\thanks{This work was supported by the Forschungsinstitut f\"{u}r Mathematik at ETH Z\"{u}rich and partially by The Leverhulme Trust}

\begin{document}

\begin{abstract} In this article, we investigate the regularity for certain elliptic systems without a $L^2$-antisymmetric structure. As applications, we prove some $\epsilon$-regularity theorems for weakly harmonic maps from the unit ball $B= B(m) \subset \mathbb{R}^m $ $(m\geq2)$ into certain pseudo-Riemannian manifolds: standard stationary Lorentzian manifolds, pseudospheres $\mathbb{S}^n_\nu \subset \mathbb{R}^{n+1}_\nu$ $(1\leq\nu \leq n)$ and pseudohyperbolic spaces $\mathbb{H}^n_\nu \subset \mathbb{R}^{n+1}_{\nu+1}$ $(0\leq\nu \leq n-1)$. Consequently, such maps are shown to be H\"{o}lder continuous (and as smooth as the regularity of the targets permits) in dimension $m=2$. In particular, we prove that any weakly harmonic map from a disc into the De-Sitter space $\mathbb{S}^n_1$ or the Anti-de-Sitter space $\mathbb{H}^n_1$ is smooth. Also, we give an alternative proof of the H\"{o}lder continuity of any weakly harmonic map from a disc into the Hyperbolic space $\mathbb{H}^n$ without using the fact that the target is nonpositively curved. Moreover, we extend the notion of generalized (weakly) harmonic maps from a disc into the standard sphere $\mathbb{S}^n$ to the case that the target is  $\mathbb{S}^n_\nu$ $(1\leq\nu \leq n)$ or $\mathbb{H}^n_\nu$ $(0\leq\nu \leq n-1)$, and obtain some $\epsilon$-regularity results for such generalized (weakly) harmonic maps.

\end{abstract}
\date{}
\date{\today}
\maketitle

\section{Introduction}
\vskip0.5cm

In the recent papers by Rivi\`{e}re \cite{Ri1} and Rivi\`{e}re-Struwe \cite{RS}, the following regularity results for elliptic systems with a $L^2$-antisymmetric structure are established:

\begin{thm}[Rivi\`{e}re \cite{Ri1} for $m=2$, Rivi\`{e}re-Struwe \cite{RS} for $m \geq 3$] Let $B= B(m) \subset \mathbb{R}^m$ $(m\geq2)$ be the unit ball. There exists $ \epsilon_m > 0$ such that for every   $\Omega \in
L^{2}(B,so(n)\otimes \wedge ^1 \mathbb{R}^{m})$ and for every weak solution $u\in W^{1,2}(B,\R^n)$ of the following elliptic system:
\bee \label{1.1} -\ {\rm div}\ \nabla u = \Omega \cdot  \nabla u \eee
satisfying
\bee \label{1.2}  \underset {B_R(x)\subset B}{\rm sup}  \left (  R^{2-m} \int_{B_R(x)} |\nabla u|^2 + |\Omega|^2 \right )^{\frac{1}{2}} < \epsilon_m,
\eee
we have that $u$ is H\"{o}lder continuous in $B$.
\end{thm}

One of the main applications of the above results is the regularity theory for harmonic map systems into closed Riemannian manifolds, where the $L^2$-antisymmetric property of the potential $\Omega$ in \eqref{1.1} relies on the fact that the target manifolds are compact and Riemannian. For classical regularity results of weakly harmonic maps, see e.g. the books by H\'{e}lein \cite{H4} and Lin-Wang \cite{LW} and references therein.

In this paper, we shall study the regularity for weakly harmonic maps from the unit ball $B= B(m) \subset \mathbb{R}^m$ $(m\geq2)$ into certain pseudo-Riemannian manifolds from different points of view. Analytically, it is interesting to know how the structure of the harmonic map system is affected when the target manifolds become pseudo-Riemannian. As we will see later, in general, the $L^2$-antisymmetric structure for harmonic map systems into closed Riemannian manifolds may not be preserved any more when the target manifolds become non-compact or non-Riemannian. Therefore, we would like to explore the extent to which the results developed by Rivi\`{e}re \cite{Ri1} and Rivi\`{e}re-Struwe \cite{RS} can be generalized to elliptic systems without a $L^2$-antisymmetric structure. Geometrically, considering the link between harmonic maps into $\mathbb{S}^4_1 \subset \mathbb{R}^5_1$ and the conformal gauss maps of Willmore surfaces in $\mathbb{S}^3$ (see Bryant \cite{Br}. See also \cite{H3, Pa, BD1, BD2}), and the regularity results for weak Willmore immersions established by Rivi\`{e}re \cite{Ri2}, we are strongly encouraged to find a method to study the regularity for weakly harmonic maps into $\mathbb{S}^4_1$ and then extend it to the cases of more general targets. Physically, it is known that harmonic maps play an important role in string theory (see e.g. \cite{De, J2}). One of the most significant results in string theory is the AdS/CFT correspondence (Anti-de-Sitter space/Conformal Field Theory correspondence) proposed in 1997 by Maldacena \cite{Ma}. In view of the recent work on minimal surfaces in Anti-de-Sitter space and its applications in theoretical physics (see e.g. Alday-Maldacena \cite{AM}), we are interested in extending the regularity theory for harmonic maps into closed Riemannian manifolds to the cases that the targets are some model spacetimes (which are non-compact and Lorentzian) considered in General Relativity (see e.g. \cite{KSHM,O}), for instance, standard stationary Lorentzian manifolds, De-Sitter space $\mathbb{S}^n_1$ (also denoted by $dS_n$) and Anti-de-Sitter space $\mathbb{H}^n_1$ (also denoted by $AdS_n$).

In the present work, we solve these problems by using the theory of integrability by compensation developed in \cite{We, M, CLMS, Fe, FS} and some conservation laws, due to the symmetries of the target manifolds considered. We point out that our results partially realize the perspectives (proposed by Rivi\`{e}re \cite{Ri2}, p.3-4) of the regularity theory for elliptic systems. For some other generalizations of the methods of Rivi\`{e}re \cite{Ri1} and Rivi\`{e}re-Struwe \cite{RS}, see Lamm-Rivi\`{e}re \cite{LR}, Struwe \cite{Stru}, Duzaar-Mingione \cite{DM} and Rivi\`{e}re \cite{Ri3}. For some other analytic aspects of harmonic maps into pseudo-Riemannian manifolds, see e.g. H\'{e}lein \cite{H5}.

First, we observe that, by slightly adapting the techniques used by Rivi\`{e}re-Struwe \cite{RS}, similar regularity
results as in Theorem 1.1 extend to certain elliptic systems with a potential a priori in $L^2$ but not necessary
antisymmetric. To see this, recall that for $1 \leq s < \infty$, the Morrey norm $|| \cdot ||_{M^s_s(B)}$ of a function
$f\in L^s_{{\rm loc}}(B)$ is
\bee || f ||_{M^s_s(B)}  =  \underset {B_R(x)\subset B}{\rm sup}  \left ( R^{s-m} \int_{B_R(x)} |f|^s \right)^{\frac{1}{s}}, \nn
\eee
then we have the following

\begin{thm}  For $m\geq2$ and for any $\Lambda > 0$, there exists $ \epsilon_{m,\Lambda} > 0$ such that for every $\Theta \in L^{2}(B,so(n)\otimes \wedge ^1 \mathbb{R}^{m})$, $\zeta  \in W^{1,2}(B, {\rm M}(n) \otimes \wedge^2 \mathbb{R}^{m})$,
 $F \in W^{1,2}\cap L^{\infty}(B, {\rm M}(n))$, $G \in W^{1,2}\cap L^{\infty}(B,{\rm M}(n))$  and $ Q \in W^{1,2}\cap L^{\infty}(B, {\rm GL}(n))$ and for every weak solution $u\in W^{1,2}(B,\R^n)$ of the following elliptic system:
\bee \label{1.3}  -\ {\rm div}\  (Q \ \nabla u)= \Theta  \cdot  Q \ \nabla u + F \ {\rm curl}\ \zeta \cdot G\  \nabla u \eee
satisfying
\bee  \label{1.4}   || \nabla u  ||_{M^2_2(B)} +   || \Theta  ||_{M^2_2(B)} +  || {\rm curl}\  \zeta  ||_{M^2_2(B)} +  || \nabla  Q  ||_{M^2_2(B)} + || \nabla F ||_{M^2_2(B)} + || \nabla  G  ||_{M^2_2(B)} <  \epsilon_{m,\Lambda}
\eee
and
\bee \label{1.5} |Q|+|Q^{-1}|+|F|+|G| \leq \Lambda, \quad {\rm a.e.\ in}\ B,
\eee
we have that $u$ is H\"{o}lder continuous in $B$ .
\end{thm}

The result in Theorem 1.2 was partially obtained by Hajlasz-Strzelecki-Zhong (\cite{HSZ}, Theorem 1.2) for the case $m=2$, $\Theta \equiv 0$, $Q \equiv I_n$ and by Schikorra (\cite{Sc}, Remark 3.4) for the case $m\geq2$, $\zeta \equiv 0$.

Note that the elliptic system \eqref{1.3} can be written as
\bee \label{1.6}  -\ {\rm div}\  (Q \ \nabla u)= \left \{\Theta  + F \ {\rm curl}\ \zeta  \ (G Q^{-1} ) \right \} \cdot  (Q \ \nabla u) \eee
or equivalently as
\bee \label{1.7}  -\ {\rm div}\  \nabla u= \left \{Q^{-1} \nabla Q + Q^{-1} \left  (\Theta  + F \ {\rm curl}\ \zeta  \ (G Q^{-1} ) \right ) Q \right \} \cdot   \nabla u. \eee
Considering $Q $ as a kind of gauge transformation,  we interpret the elliptic system \eqref{1.7} as follows: its potential $$Q^{-1} \nabla Q + Q^{-1} \left  (\Theta  + F \ {\rm curl}\ \zeta  \ (G Q^{-1} ) \right ) Q  $$  is gauge equivalent to a new one $$  \Theta  + F \ {\rm curl}\ \zeta  \ (G Q^{-1} )$$ which can be decomposed into an antisymmetric part $\Theta$  and an almost divergence free part $F \ {\rm curl}\ \zeta  \ (G Q^{-1} )$.

As an application of Theorem 1.2, we shall study the regularity for weakly harmonic maps into standard stationary
Lorentzian manifolds. A standard stationary Lorentzian manifold (see e.g. \cite{KSHM,O}) is a product manifold $\mathbb{R} \times M$ equipped with a metric
\bee \label{1.8} g = -\ (\beta\circ \pi_M ) \left (\pi^*_\mathbb{R}dt + \pi^*_M \omega \right ) \otimes \left (\pi^*_\mathbb{R}dt + \pi^*_M \omega \right ) +  \pi^*_M g_M, \eee
where $(\mathbb{R}, dt^2)$ is the 1-dimensional Euclidean space, $(M, g_M)$ is a closed Riemannian manifold of class  $C^3$, $\beta $ is a  positive $C^2$ function on $M$, $\omega$
 is a $C^2$ 1-form on $M$, $\pi_\mathbb{R}$ and $\pi_M $ are the natural projections on $\mathbb{R}$ and $M$, respectively. For simplicity of notations, we shall write the metric \eqref{1.8}  as
\bee \label{1.9} g = -\  \beta(dt + \omega)^2 + g_M. \eee
By Nash's embedding theorem, we embed $(M,g_M)$ isometrically into some Euclidean space $\mathbb{R}^n$. Then, there exists a tubular neighborhood $V_{\delta}M$ of radius $\delta>0 $ of $M$ in $\mathbb{R}^n$ and a $C^2$ projection map $\Pi$ from  $V_{\delta}M$ to $M$ (see H\'{e}lein's book \cite{H4}, Chapter 1). Moreover, we pull back $\b$ and $\omega$ via the projection $\Pi$ and obtain  $\Pi^* \b \in C^2(V_{\delta}M, (0,\infty))$ and $\Pi^*\omega \in C^2(\Omega^1(V_{\delta}M)) $, respectively. For simplicity, we shall still denote $ \Pi^* \b$ and $\Pi^* \omega$ by $\b$ and $\omega$, respectively. Write $\omega = \sum_{i=1}^n \omega_i(y)dy^i$, $y=(y^1,...,y^n)\in V_{\delta}M \subset \mathbb{R}^n$,where $\omega_i \in C^2(V_{\delta}M) $.

To study the regularity for weakly harmonic maps into $(\mathbb{R}\times M,g)$, we consider the space
\bee W^{1,2}(B,\mathbb{R}\times M ):= \left \{ \  (t,u) \in W^{1,2}(B, \mathbb{R} ) \times  W^{1,2}(B,\mathbb{R}^{n} ) \   |  \  u(x) \in M \ { \rm a.e.}\  x \in B \  \right \}
\eee
For a map $(t,u)\in W^{1,2}(B, \mathbb{R}\times M )$, we define the following Lagrangian:
\bee \label{1.11} E(t,u)= -   \frac{1}{2}\int_B \beta(u) \left |\nabla t  +  \omega_i(u)\nabla u^i \right |^2 +  \frac{1}{2}\int_B   |\nabla u |^2
\eee

\begin{defi} A map $(t,u) \in W^{1,2}(B,\mathbb{R}\times M)$ is called a weakly harmonic map from $B$ into $(\mathbb{R}\times M,g)$, if it is a critical point of the Lagrangian functional  \eqref{1.11}.
\end{defi}

The Euler-Lagrange equation (see Section 3) for a weakly harmonic map $(t,u)\in W^{1,2}(B, \mathbb{R}\times M )$ from $B$ into $(\mathbb{R}\times M,g)$ is an elliptic system of the form \eqref{1.3}, which can be geometrically
interpreted as follows: the antisymmetric term $\Theta$ corresponds to the Riemannian structure of the closed
spacelike hypersurfaces $\{t\} \times M$ and the divergence free term ${\rm curl}\ \zeta$ corresponds to the following conservation law
\bee \label{1.12}
 {\rm div} \  \left\{ \ \beta(u) (\nabla t  +  \omega_i(u) \nabla u^i  ) \  \right \} = 0,  \quad {\rm in}\ \mathcal{D}'(B),
\eee
due to the symmetry of the target generated by the timelike Killing vector field $\partial_t$. Applying Theorem 1.2, we have the following $\epsilon$-regularity result:

\begin{thm} For $m\geq 2$, there exists $\epsilon_{m}>0 $ depending on $(\mathbb{R}\times M,g)$ such that any weakly harmonic map $(t,u) \in W^{1,2}(B,\mathbb{R}\times M)$ from $B$ into $(\mathbb{R}\times M,g)$ satisfying
\bee \label{1.13}  || \nabla t  ||_{M^2_2(B)} + || \nabla u  ||_{M^2_2(B)} < \epsilon_{m},
\eee
is  H\"{o}lder continuous (and as smooth as the regularity of the target permits) in $B $.
\end{thm}

In dimension $m=2$, we notice that the Morrey norm $|| \cdot ||_{M^2_2}$ reduces to the norm $|| \cdot ||_{L^2}$. Therefore, by conformal invariance and rescaling in the domain, we obtain the following regularity result:

\begin{thm}  For $m=2$, any weakly harmonic map $(t,u) \in W^{1,2}(B,\mathbb{R}\times M)$ from $B$ into $(\mathbb{R}\times M,g)$ is  H\"{o}lder continuous (aand as smooth as the regularity of the target permits) in $B$.
\end{thm}

In Theorem 1.4, if the target $(\mathbb{R}\times M,g)$ is a standard static Lorentzian manifold (see e.g. \cite{KSHM,O}), namely, the 1-form $\omega$ in the metric $g$ (see \eqref{1.9}) vanishes identically, then the corresponding regularity result was proved by Isobe \cite{Is} (using H\'{e}lein's method of moving frame \cite{H4}).

\vskip 0.2cm

Next, we shall consider, in a certain sense, elliptic systems of the form \eqref{1.1} with the potential $\Omega$ a priori only in $ L^p$ for some $1<p<2$. Note that, if $\Omega$ is not in $L^2$, then the right hand side of \eqref{1.1} is not in $L^1$ and thus the equation makes no sense any more (not even in the distribution sense!). However, we observe that, if in addition, $\Omega$ is divergence free,
namely,
\bee \label{1.14} {\rm div }\ \Omega  = 0, \quad {\rm in}\ \mathcal{D}'(B),
\eee
then the equation \eqref{1.1} can be written in the following form:
\bee  \label{1.15}  -\ {\rm div}\  \left ( \nabla u + \Omega  \  u  \right ) = 0, \quad {\rm in}\ \mathcal{D}'(B),
\eee
This new form \eqref{1.15} has the advantage that it is still meaningful in the distribution sense if $\Omega$ is a priori only in $L^p$ for some $1<p<2$. Moreover, under the further assumption that the Morrey norms $|| \nabla u ||_{M^p_p(B)}$ and  $|| \Omega ||_{M^p_p(B)}$ are sufficiently small, the H\"{o}lder continuity of the weak solution $u$ holds.

\begin{thm} For $m\geq2$ and for any  $1<p<\frac{m}{m-1} $,
 there exists $\epsilon_{m,p} > 0$ such that for any $\Omega \in L^p(B, {\rm M}(n)\otimes \wedge ^1 \mathbb{R}^{m})$
satisfying \eqref{1.14} and for any weak solution  $u\in W^{1,2}( B ,\mathbb{R}^n)$ of the elliptic system \eqref{1.15} satisfying
\bee  \label{1.16} ||  \nabla u || _{ M^p_p(B) } + ||  \Omega || _{ M^p_p(B) }  <  \epsilon_{m,p},
\eee
we have that $u$ is H\"{o}lder continuous in $B$.
\end{thm}

As applications of Theorem 1.5, we shall study the regularity for weakly harmonic maps into pseudospheres and pseudohyperbolic spaces. For this purpose, we recall some facts about these target spaces and refer to O'Neill's book \cite{O} for more details.

Let $n\in \mathbb{N}$ and let $\nu \in \mathbb{N}$ satisfy $0 \leq \nu \leq n$. Denote
\bee \label{1.17} \mathcal {E} = (\varepsilon_{ij}) : =  \left ( \ba{ccc}
            -I_{\nu} & 0   \\
            0 &  I_{n+1 -\nu}  \ea \right ).
\eee
The pseudo-Euclidean space $\mathbb{R}^{n+1}_\nu$
of signature $(\nu, n+1-\nu)$ is the space $\mathbb{R}^{n+1}$ equipped with a metric
\bee  \langle v, w \rangle_{\mathbb{R}^{n+1}_\nu} := v^T\mathcal{E}\ w =  -\ \left ( v^1w^1+ ... + v^\nu w^\nu \right )+  \left (v^{\nu+1} w^{\nu+1} +... + v^{n+1}w^{n+1} \right ),
\nn\eee
for all $v=(v^1, ... , v^{n+1})^T \in \mathbb{R}^{n+1}$ and $w=(w^1, ... ,w^{n+1})^T \in \mathbb{R}^{n+1}$. The pseudoshpere  $\mathbb{S}^n_\nu$ in $\mathbb{R}^{n+1}_\nu$ is defined as
\bee \mathbb{S}^n_\nu &: =&  \left \{ \ y \in \mathbb{R}^{n+1}_\nu \mid  \langle y, y \rangle_{\mathbb{R}^{n+1}_\nu} = y^T \mathcal {E}\ y= 1  \  \right \}
\eee
with the induced metric. In particular, $\mathbb{S}^n_0 \subset \mathbb{R}^{n+1}_0$ is the standard sphere $\mathbb{S}^n \subset \mathbb{R}^{n+1}$
 and $\mathbb{S}^n_1 \subset \mathbb{R}^{n+1}_1$ is the De-Sitter space $dS_n$ in General Relativity. The linear isometries of $\mathbb{R}^{n+1}_\nu$ form the group
\bee \label{1.19} {\rm O}(\nu, n+1-\nu) = \left \{ \ P\in {\rm GL}(n+1) \mid  P^T = \mathcal {E}  P^{-1} \mathcal {E} \ \right \}.
\eee
Denote by ${\rm SO}^+(\nu, n+1-\nu)$ the identity component of ${\rm O}(\nu, n+1-\nu)$. The lie algebra of
${\rm SO}^+(\nu, n+1 -\nu)$ is
\bee so(\nu, n+1 -\nu) = \left \{  \ A \in {\rm GL}(n+1)  \mid  A^T = - \mathcal {E}\  A\ \mathcal {E}  \  \right \}.
\eee

Using the isometric embedding $\mathbb{S}^n_\nu \subset \mathbb{R}^{n+1}_\nu$, we set
\bee \label{1.21} W^{1,2}(B,\mathbb{S}^n_\nu ):= \left  \{\ u = (u^1,u^2, ... , u^{n+1} )^T \in W^{1,2}(B,\mathbb{R}^{n+1}_\nu ) \mid u^T \mathcal {E}\ u =1\  { \rm a.e.\  in}\  B \  \right \}.
\eee
For a map $u\in W^{1,2}(B,\mathbb{S}^n_\nu )$, we define the following Lagrangian:
\bee \label{1.22}  E(u) := \frac{1}{2}\int_B   (\nabla u )^T  \mathcal{E}\ \nabla u = -  \frac{1}{2}\int_B  \left ( |\nabla u^1|^2 + ... +  |\nabla u^\nu|^2 \right ) +  \frac{1}{2}\int_B \left ( |\nabla u^{\nu+1}|^2 +... + |\nabla u^{n+1}|^2 \right )
\eee

\begin{defi}  A map $u\in W^{1,2}(B,\mathbb{S}^n_\nu )$ is called a weakly harmonic map from $B$ into $\mathbb{S}^n_\nu$, if it is a critical point of the Lagrangian functional \eqref{1.22}.
\end{defi}

Denote
\bee \label{1.23} \mathcal {F} = (\varsigma_{ij}) := \left ( \ba{ccc}
            -I_{\nu+1 } & 0   \\
            0 &  I_{n -\nu}  \ea \right ).
\eee
The pseudohyperbolic space  $\mathbb{H}^n_\nu$ in  $\mathbb{R}^{n+1}_{\nu+1}$ is defined as
\bee \mathbb{H}^n_\nu &: =& \left \{ y \in \mathbb{R}^{n+1}_{\nu+1} \mid  \langle y, y \rangle_{\mathbb{R}^{n+1}_{\nu+1}} = y^T \mathcal {F} y  = - 1  \right \}
\eee
with the induced metric. In particular, $\mathbb{H}^n_0  \subset \mathbb{R}^{n+1}_1$ is a hyperboloid containing two copies of the Hyperbolic space $\mathbb{H}^n $ and $\mathbb{H}^n_1  \subset \mathbb{R}^{n+1}_2$ is the Anti-de-Sitter space $AdS_n$ in General Relativity.

Using the isometric embedding $\mathbb{H}^n_\nu \subset \mathbb{R}^{n+1}_{\nu+1}$, we set
\bee \label{1.25} W^{1,2}(B,\mathbb{H}^n_\nu ):= \left  \{\ u = (u^1,u^2, ... , u^{n+1} )^T \in W^{1,2}(B,\mathbb{R}^{n+1}_{\nu+1}) \mid u^T \mathcal {F} u = - 1\  { \rm a.e.\  in}\  B \ \right  \}.
\eee
For a map $u\in W^{1,2}(B,\mathbb{H}^n_\nu )$, we define the following Lagrangian:
\bee \label{1.26} E(u) := \frac{1}{2}\int_B   (\nabla u )^T  \mathcal{F}\ \nabla u = -  \frac{1}{2}\int_B  \left ( |\nabla u^1|^2 + ... +  |\nabla u^{\nu+1}|^2 \right ) +  \frac{1}{2}\int_B \left ( |\nabla u^{\nu+2}|^2 +... + |\nabla u^{n+1}|^2 \right )
\eee

\begin{defi} A map $u\in W^{1,2}(B,\mathbb{H}^n_{\nu})$ is called a weakly harmonic map from $B$ into $\mathbb{H}^n_{\nu} $, if it is a critical point of the Lagrangian functional \eqref{1.26}.
\end{defi}

Notice that the following anti-isometry (see O'Neill's book \cite{O})
\bee \sigma: \hskip1.5cm  \mathbb{R}^{n+1}_\nu &\rightarrow& \mathbb{R}^{n+1}_{n-\nu+1}  \nn \\
                         (y_1,...,y_{n+1})  &\mapsto& (y_{\nu+1},...,y_{n+1}, y_1,...,y_\nu)   \nn
\eee
induces an anti-isometry from $\mathbb{S}^n_\nu$ to $\mathbb{H}^n_{n-\nu}$.
In the sequel, we shall only consider the cases of $\mathbb{S}^n_\nu$ $(0\leq \nu \leq n)$.

To proceed, we recall that a weakly harmonic map $u\in W^{1,2}(B,\mathbb{S}^n)$ satisfies the following conservation laws (due to Shatah \cite{Sha} and Chen \cite{Che}. See also Rubinstein-Sternberg-Keller \cite{RSK} and H\'{e}lein's book \cite{H4}):
\bee  \label{1.27} {\rm div } \left (u^i \nabla u^j - u^j \nabla u^i \right ) = 0,\quad {\rm in}\ \mathcal{D}'(B), \quad \forall i,j=1,2,...,n+1,
\eee
which can be interpreted by Noether theorem, using the symmetries of $\mathbb{S}^n$. Note that the pseudospheres $\mathbb{S}^n_\nu$ $(1\leq \nu \leq n)$ have isometry groups $ {\rm O}(\nu, n+1-\nu)$ and hence they are all maximally symmetric. With the help of the symmetric properties, we are able to extend the conservation laws \eqref{1.27} to weakly harmonic maps into these more general targets.

\begin{pro} For $ m \geq 2$. Let $u\in W^{1,2}(B,\mathbb{S}^n_\nu )$ $(1\leq \nu \leq n)$ be a weakly harmonic map. Then the conservation laws \eqref{1.27} hold.
\end{pro}

For a weakly harmonic map $u\in W^{1,2}(B,\mathbb{S}^n_\nu )$ $(0\leq \nu \leq n)$, we define the following matrix valued vector field
\bee \label{1.28} \Theta = \left (\Theta^{ij} \right ) :=  \left (u^i \nabla u^j - u^j \nabla u^i \right ), \quad  i,j=1,2,...,n+1. \eee

In the case of a compact target $\mathbb{S}^n$, $\Theta$ is in  $L^2(B, {\rm M}(n)\otimes \wedge ^1 \mathbb{R}^{m})$ and $u$ weakly solves the following elliptic system (see H\'{e}lein \cite{H1, H4})
\bee   - \ {\rm div}\  \nabla u =  \Theta \cdot  \nabla u. \nn
\eee
Since $\Theta$ is divergence free (due to the conservation laws \eqref{1.27}), the continuity of $u$ in dimension $m=2$ follows immediately from Wente's lemma \cite{We}.

However, in the case of a non-compact target $\mathbb{S}^n_\nu$ $(1\leq \nu \leq n)$, $\Theta$  is only in  $L^p(B, {\rm M}(n)\otimes \wedge ^1 \mathbb{R}^{m})$ for any $1<p<2$. Proposition 1.1 indicates that $\Theta$ is still divergence free. In what follows, we show that $u$ is a weak solution of an elliptic system of the form \eqref{1.15} with its potential satisfying \eqref{1.14}. Moreover, by making use of the conservation laws \eqref{1.27}, we are able to estimate $||\Theta||_{M^p_p(B_{1/2})}$ by $|| \nabla u ||_{M^p_p(B)}$, where $B_{1/2}= B_{1/2}(m) \subset \mathbb{R}^m$ $(m\geq2)$ is the ball centered at 0 and of radius $1/2$.

\begin{pro} For $ m \geq 2$. Let $u\in W^{1,2}(B,\mathbb{S}^n_\nu )$ $(1\leq \nu \leq n)$. Then
\bee \label{1.29}   \nabla u + \Theta \ \mathcal {E}\  u  = 0,  \quad {\rm a.e. \ in}\ B.
\eee
where $\Theta$ is defined as in \eqref{1.28}. Consequently, we have
\bee \label{1.30}   - \ {\rm div} \left ( \nabla u + \Theta \ \mathcal {E}\  u \right ) = 0  \quad {\rm in}\ \mathcal{D}'(B).
\eee
 Furthermore, suppose that $u$ is weakly harmonic and for any fixed $1<p< \frac{m}{m-1}$ there holds $||\nabla u ||_{ M^p_p (B)} <  \infty$, then we have the following estimate:
\bee \label{1.31}
|| \Theta ||_{M^p_p(B_{1/2})} \leq C \ || \nabla u ||^2_{ M^p_p (B)}.
\eee
\end{pro}

Since $\mathcal {E}$  is a constant matrix, applying Theorem 1.5 with $\Omega =\Theta \ \mathcal {E}$ and using a rescaling of the domain gives the following $\epsilon$-regularity result:
\begin{thm} For $ m \geq 2$ and for any $1<p< \frac{m}{m-1} $, there exists $\epsilon_{m,p}>0 $ such that
 any weakly harmonic map $u\in W^{1,2}(B,\mathbb{S}^n_\nu )$ $(1\leq \nu \leq n)$ satisfying
\bee
 || \nabla u ||^2_{ M^p_p (B)} < \epsilon_{m,p}
\eee
is H\"{o}lder continuous (and hence smooth) in $B$.
\end{thm}

In dimension $m=2$,  a straightforward calculation gives that $|| \nabla u ||_{ M^p_p (B)} \leq  || \nabla u ||_{ L^2 (B)}$ for any $1<p<2$. Therefore, by conformal invariance, we have
\begin{thm} For $m=2$, any weakly harmonic map $u\in W^{1,2}(B,\mathbb{S}^n_\nu )$ $(1\leq \nu \leq n)$ is H\"{o}lder continuous (and hence smooth) in $B$.
\end{thm}

In particular, we prove that any weakly harmonic map from a disc into the De-Sitter space $\mathbb{S}^n_1$ or the Anti-de-Sitter space $\mathbb{H}^n_1 \cong \mathbb{S}^n_{n-1}$ is smooth. Also, we give an alternative proof of the H\"{o}lder continuity of weakly harmonic maps from a disc into the Hyperbolic space $\mathbb{H}^n$ (one component of $\mathbb{H}^n_0 \cong \mathbb{S}^n_{n}$) without using the fact that the target has non-positive sectional curvature (for a proof using the curvature property, we refer to Jost's book \cite{J1}). We expect that the results in Theorem 1.6 and Theorem 1.7 can be extended in the same spirit of H\'{e}lein's setting in \cite{H2} to certain homogeneous pseudo-Riemannian manifolds.

\vskip 0.2cm

Furthermore, we observe that the methods used in the proofs of Proposition 1.2 and Theorem 1.5 can be applied to study the $\epsilon$-regularity of maps in the spaces of distributions of lower regularity. This motivates us to extend the notion of generalized (weakly) harmonic maps from a disc into the standard sphere $\mathbb{S}^n$ (introduced by Almeida \cite{Al}) to the cases that the targets are pseudospheres $\mathbb{S}^n_\nu$ $(1\leq \nu \leq n)$ (see Section 5). To see this, we recall the notion of generalized (weakly) harmonic maps into $\mathbb{S}^n$.

\begin{defi}[Almeida \cite{Al}]  For $m=2$, a map $u\in W^{1,1}(B,\mathbb{S}^n )$ is called a generalized (weakly) harmonic map if \eqref{1.27} hold.
\end{defi}

Generalized (weakly) harmonic maps into $\mathbb{S}^n$ might be not continuous. However, there are some $\epsilon$-regularity results for such maps. Almeida \cite{Al} showed that any generalized harmonic map $u\in W^{1,1}(B,\mathbb{S}^n )$ with $||\nabla u||_{L^{(2,\infty)}}$ small is smooth (an alternative proof was given by Ge \cite{Ge}). Moser \cite{Mo} proved that any generalized harmonic map $u\in W^{1,p}_{{\rm loc}}(B,\mathbb{S}^n )$ with $p\in (1,2)$ is smooth if $p$ is sufficiently close to 2 and $||u||_{{\rm BMO}}$ is small. Strzelecki \cite{Strz2} showed that any generalized harmonic map $u\in W^{1,p}_{{\rm loc}}(B,\mathbb{S}^n )$ with $p\in (1,2)$ is smooth provided that $|| u||_{{\rm BMO}}$ is small.

To extend the notion of generalized (weakly) harmonic maps into the pseudospheres $\mathbb{S}^n_\nu$ $(1\leq \nu \leq n)$, we observe that a $W^{1,1}$ map from a disc into any of these non-compact targets is not a priori in $L^{\infty}$ and hence the conservation laws \eqref{1.27} make no sense for such a map. Therefore, we need to require that the map $u$ belongs to the sobolev space $ W^{1,\frac{4}{3}}$ so that $$u^i \nabla u^j - u^j \nabla u^i \in L^1_{{\rm loc}}(B),  \quad \forall i,j=1,2,...,n+1.$$ and hence the conservation laws \eqref{1.27} become meaningful.

\begin{defi} For $m=2$, a map $u\in W^{1,\frac{4}{3}}(B, \mathbb{S}^n_\nu)$ $(1\leq \nu \leq n )$ is called a generalized (weakly) harmonic map if \eqref{1.27} hold.
\end{defi}

Analogously to Theorem 1.6, we have the following $\epsilon$-regularity result.

\begin{thm} For $m=2$ and for any $\frac{4}{3}<p< 2$, there exists $\epsilon_{p}>0 $ such that any generalized (weakly) harmonic map $u\in W^{1,\frac{4}{3} }(B,\mathbb{S}^n_\nu)$ $(1\leq \nu \leq n )$ satisfying
\bee
 || \nabla u ||^2_{ M^p_p (B)} < \epsilon_{p}
\eee
is H\"{o}lder continuous (and hence smooth) in $B$.
\end{thm}

Finally, we study the regularity for an elliptic system of the form \eqref{1.1} with  $\Omega\in  L^{2}(B,so(1,1)\otimes \wedge ^1 \mathbb{R}^{2})$ in dimension $m=2$ and show by constructing an example that weak solutions in $W^{1,2}$ to such an elliptic system might be not in $L^{\infty}$.

The paper is organized as follows. In Section 2, we prove Theorem 1.2 and Theorem 1.5. In Section 3, we apply Theorem 1.2 to prove the $\epsilon$-regularity (Theorem 1.3) of weakly harmonic maps into standard stationary Lorentzian manifolds. In Section 4, we first show Proposition 1.1 and Proposition 1.2. Then we prove the regularity results (Theorem 1.6 and Theorem 1.7) for weakly harmonic maps into pseudospheres. In Section 5, the $\epsilon$-regularity result (Theorem 1.8) for generalized (weakly) harmonic maps from a disc into pseudospheres is proved. In Section 6, we study an elliptic system with a $L^2$-$so(1,1)$ structure in dimension $m=2$.

\vskip 0.2cm

{\noindent}{\bf Notation:} For a 2-vector field $\xi=\xi_{ij}\partial_{x_i} \wedge \partial_{x_j}$, ${\rm curl}\ \xi$ denotes the vector field $\left(\sum_i \partial_{x_i} \xi_{ij} \right)\partial_{x_j}$ and $d \xi$ denotes the 3-vector field $\left (\partial_{x_k} \xi_{ij} \right ) \ \partial_{x_k} \wedge \partial_{x_i} \wedge \partial_{x_j}$. A constant $C$ may depend on $m$, $n$ and $p$.

\vskip 0.2cm

{\noindent}{\bf Acknowledgements} The author would like to thank Professor Tristan Rivi\`{e}re for many useful  discussions, consistent support and encouragement. Thanks also to Professor J\"{u}rgen Jost and Professor Yuxin Ge for helpful conversations.

\vskip1.5cm

\section{Proofs of Theorem 1.2 and Theorem 1.5}
\vskip0.5cm

In this section, we will prove Theorem 1.2 and Theorem 1.5.

First, combining the div-curl inequality by Coifman-Lions-Meyer-Semmes \cite{CLMS} (see M\"{u}ller \cite{M} for an earlier contribution), the Hardy-BMO duality by Fefferman \cite{Fe} (see also Fefferman-Stein \cite{FS} and Stein \cite{Ste}) and the observation that the Morrey spaces $M^s_s(\mathbb{R}^m)$ ($ 1 \leq s < \infty$) are contained in the space BMO($\mathbb{R}^m$) (due to Evans \cite{Ev}), we give the following lemma (see Proposition III.2 in Bethuel \cite{Be}, Lemma 3.1 in Schikorra \cite{Sc} and Strzelecki \cite{Strz1} p.234-235. See also Chanillo \cite{Ch} and Chanillo-Li \cite{CL}).

\begin{lem}  For $m \geq 2$, $1 \leq s < \infty$ and $ 1< p < \infty$.
Let $1< q < \infty$ satisfy $\frac{1}{p} + \frac{1}{q}=1$. For any ball
 $B_R(x)\subset \mathbb{R}^m$, $f \in W^{1,p}(B_{R}(x))$, $g \in W^{1,q}(B_{R}(x), \wedge^2\mathbb{R}^m)$ satisfying
\bee f|_{\partial B_R(x)} = 0  \quad \ {\rm or} \quad  g|_{\partial B_R(x)} = 0
\eee
and $h \in W^{1,s}(B_{2R}(x)) $ satisfying
\bee || \nabla h||_{M^s_s(B_{2R}(x))} < \infty,
\eee
there holds:
\bee \int_{B_{R}(x)} \left ( \nabla f \cdot \ {\rm curl}\ g  \right ) \ h \leq  C \ ||\nabla f ||_{L^p(B_R(x))} \ ||{\rm curl}\ g ||_{L^q(B_R(x))} \ || \nabla h||_{M^s_s(B_{2R}(x))},
\eee
where $C=C_{m,s,p}>0$ is a uniform constant independent of $R>0$.
\end{lem}

Next, with the help of the above lemma, we follow the approach used by Rivi\`{e}re-Struwe \cite{RS} to prove Theorem 1.2 and Theorem 1.5.

\vskip0.2cm

\noindent{\bf Proof of Theorem 1.2:} Fix $m\geq 2$ and $\Lambda >0$. Choose $\epsilon_{m,\Lambda}>0$ sufficiently small, then by assumption \eqref{1.4} and the existence of Coulomb gauge (due to Rivi\`{e}re \cite{Ri1} for $m=2$ and Rivi\`{e}re-Struwe \cite{RS} for $m \geq 3$), we conclude that there are $P\in W^{1,2}(B, {\rm SO}(n))$ and $\xi \in
W^{1,2}_0(B,so(n) \otimes \wedge^2\mathbb{R}^m)$  with $d \xi =0$ such that
\bee \label{2.4} P^{-1}\nabla P+P^{-1}\Theta P &=& {\rm curl}\  \xi \quad {\rm in}\ B.
\eee
and the following estimate holds
\bee \label{2.5} || \nabla P ||_{M^2_2(B)} + \parallel\nabla \xi
\parallel_{M^2_2(B)} \leq C
\parallel \Theta \parallel_{M^2_2(B)} \leq C \epsilon_{m,\Lambda}. \eee

Using \eqref{2.4}, we rewrite the system \eqref{1.3} as
\bee   -\ {\rm div}\  \left (P^{-1}Q \ \nabla u \right ) &=& \left ( P^{-1}\nabla P + P^{-1}\Theta P \right ) \cdot P^{-1}Q \ \nabla u   + P^{-1}F \ {\rm curl}\ \zeta \ G  \cdot  \nabla u \nn\\
&=&    {\rm curl}\ \xi  \cdot P^{-1}Q \ \nabla u   + P^{-1}F \ {\rm curl}\ \zeta \ G  \cdot  \nabla u
 \eee
Write $P^{-1}= (P_{ij}), \Theta = (\Theta^{ij})$, $\zeta = (\zeta^{ij})$,
 $F = (F^{ij})$, $G=(G^{ij})$  and $ Q=(Q^{ij}) $. Then the above equation can be written as
\bee \label{2.7}  -\ {\rm div}\  ( \sum_{j,k}  P_{ij} Q^{jk} \ \nabla u^k) &=& \sum_{j,k,l} {\rm curl}\ \xi^{ij} \cdot P_{jk} Q^{kl} \ \nabla u^l + \sum_{j,k,l,r} P_{ij}F^{jk}  {\rm curl}\ \zeta^{kl} \cdot G^{lr}\  \nabla u^r \nn\\
&=& \sum_{j,k,l}  P_{jk} Q^{kl} \ {\rm curl}\  \xi^{ij} \cdot  \nabla u^l + \sum_{j,k,l,r} P_{ij}F^{jk} G^{lr}  {\rm curl}\ \zeta^{kl} \cdot  \nabla u^r
\eee
Since $P\in W^{1,2}(B, {\rm SO}(n))$,$F \in W^{1,2}\cap L^{\infty}(B, {\rm M}(n))$, $G \in W^{1,2}\cap L^{\infty}(B,{\rm M}(n))$  and $ Q \in W^{1,2}\cap L^{\infty}(B, {\rm GL}(n))$,
we have $P^{-1}Q   \in W^{1,2}\cap L^{\infty}(B,{\rm GL}(n))$,  $P_{ij}F^{jk} G^{lr} \in W^{1,2}\cap L^{\infty}(B)$. Using the assumption \eqref{1.5}, one can verify that
\bee   \label{2.8}  || \nabla (P^{-1}Q ) ||_{M^2_2(B)} + \sum_{i,k,l,r} || \nabla (P_{ij}F^{jk} G^{lr}) ||_{M^2_2(B)}  \leq  C(\Lambda) \left (  || \nabla P ||_{M^2_2(B)} +  || \nabla Q ||_{M^2_2(B)} + || \nabla F ||_{M^2_2(B)} + || \nabla G ||_{M^2_2(B)} \right ).
\eee
Here and in the sequel, $C(\Lambda)>0$ is a constant also depending on $\Lambda$.

Combining \eqref{2.5}, \eqref{2.8} and assumption \eqref{1.4}, we get
\bee  \label{2.9} || \nabla u ||_{M^2_2(B)} + || \nabla (P^{-1}Q) ||_{M^2_2(B)} + \sum_{i,j,k,l,r} || \nabla (P_{ij}F^{jk} G^{lr}) ||_{M^2_2(B)}
+   || \nabla \xi ||_{M^2_2(B)} + || {\rm curl}\  \zeta  ||_{M^2_2(B)} \leq  C(\Lambda) \epsilon_{m,\Lambda}.
\eee

On the other hand, since $P^{-1}$ takes values in ${\rm SO} (n)$, it follows from assumption \eqref{1.5} that
\bee  \label{2.10}  C(\Lambda)^{-1}\   |\nabla u|\leq |P^{-1}Q \ \nabla u| =  |Q \ \nabla u| \leq C(\Lambda) \  |\nabla u|, \quad {\rm a.e.\ in}\ B.
\eee

Similarly to the approach by Rivi\`{e}re-Struwe (\cite{RS}, Section 3, Proof of Theorem 1.1, p.459-460. See also Schikorra \cite{Sc}, p.510-511), we apply Hodge decomposition (see \cite{IM}) to $ P^{-1}Q \ \nabla u $, use \eqref{2.7}, \eqref{2.9}, \eqref{2.10}, Lemma 2.1, and take $\epsilon_{m,\Lambda}>0$ sufficiently small to get the Morrey type estimates for $\nabla u$. Finally, we apply an iteration argument as in \cite{Gi} to obtain the H\"{o}lder continuity of $u$ in $B$.   \eop

\vskip0.2cm

\noindent{\bf Proof of Theorem 1.5:} Fix any $1<p<\frac{m}{m-1}$. Since $ {\rm div }\ \Omega  = 0 $, by Hodge decomposition, there exists $\xi \in W^{1,p}(B,  {\rm M}(n) \otimes \wedge^2 \mathbb{R}^{m}  )$ such that
\bee \label{2.11}  \Omega  =  {\rm curl}\ \xi.
\eee

Let $B_{2R}(x_0) \subset B$ and let $w \in W^{1,2}(B_R(x_0),\mathbb{R}^n) $ be solving
\bee     \label{2.12}
\left\{
\ba{rcll}
-\  {\rm div}\  \nabla w  &=&   0, & \qquad  \text{ in } B_R(x_0)  \\
w &=& u,       &  \qquad    \text{ on }    \partial B_R(x_0)
\ea
\right.
\eee
Then $v:=u-w \in W^{1,2}(B_R(x_0),\mathbb{R}^n) $ solves
\bee \label{2.13}
\left
\{    \ba{rcll}
  - \ {\rm div}\ ( \nabla v  +  \Omega\  u ) &=& 0 ,   & \qquad  \text{ in } B_R(x_0)    \\
               v &=& 0,        &  \qquad    \text{ on }   \partial B_R(x_0) \ea
\right.
\eee

Let $q = \frac{p}{p-1} > m $ be the conjugate exponent of $p$. For any $\varphi\in W^{1,q}_0(B_R(x_0))$ with $||\varphi||_{ W^{1,q}(B_R(x_0))} \leq 1$.
Using assumption \eqref{1.16}, Lemma 2.1,  \eqref{2.11} and \eqref{2.13}, we estimate for each $i$,
\bee \label{2.14}  \  \int_{B_R(x_0)}  \nabla v^i \cdot \nabla \varphi  &=& -  \int_{B_R(x_0)}  (\Omega^{ij}\ u ^j) \cdot  \nabla \varphi  \nn\\
&=& -  \int_{B_R(x_0)} ( {\rm curl}\  \xi^{ij}  \cdot  \nabla \varphi) \ u^j   \nn\\
&\leq& C ||  {\rm curl}\ \xi^{ij} ||_{ L^p(B_R(x_0))}  \ ||  \nabla \varphi|| _{ L^q(B_R(x_0))} \    ||  \nabla u || _{M^p_p(B_{2R}(x_0)) }   \nn\\
&\leq&  C  ||  \Omega^{ij} ||_{ L^p(B_R(x_0))}  \   ||  \nabla u || _{M^p_p(B_{2R}(x_0)) } \nn\\
&\leq&  C  R^{\frac{m}{p}-1}  ||  \Omega^{ij} || _{ M^p_p(B) } \  ||  \nabla u || _{M^p_p(B_{2R}(x_0)) }      \nn\\
&\leq&  C  R^{\frac{m}{p}-1} \epsilon_{m,p}  \ ||  \nabla u || _{M^p_p(B_{2R}(x_0)) }.
\eee
Since $v |_{\partial B_R(x_0)} =0 $, by duality (similarly to Rivi\`{e}re-Struwe \cite{RS}) there holds:
\bee  \label{2.15}   || \nabla v||_{L^p(B_R(x_0))}  &\leq&   C\  \underset{\varphi\in W^{1,q}_0(B_R(x_0)),\  ||\varphi||_{ W^{1,q}} \leq 1} {\rm sup } \int_{B_R(x_0)}  \nabla v \cdot \nabla \varphi.
\eee
Combining \eqref{2.14} and \eqref{2.15} gives
\bee    \label{2.16}   || \nabla v||_{L^p(B_R(x_0))}  &\leq&   C\  R^{\frac{m}{p}-1} \epsilon_{m,p} \ ||  \nabla u || _{M^p_p(B_{2R}(x_0)) }.
\eee

Next, we see from \eqref{2.12} that $w $ is harmonic in $B_R(x_0) $ and hence $\nabla w$ is also harmonic in $B_R(x_0) $. By Campanato estimates for harmonic functions (see \cite{Gi}), we have that for any $r<R$ the following holds:
\bee   \label{2.17} \int_{B_r(x_0)} |\nabla w|^p \leq    C\  \left (\frac{r}{R} \right )^m  \int_{B_R(x_0)} |\nabla w|^p.
\eee

Using that fact that $u=v+w$ and combining \eqref{2.16}, \eqref{2.17}, we estimate
\bee \label{2.18}   \int_{B_r(x_0)} |\nabla u|^p  &\leq&   C  \int_{B_r(x_0)} |\nabla w|^p   +   C  \int_{B_r(x_0)} |\nabla v|^p  \nn\\
&\leq&  C \left (\frac{r}{R} \right)^m  \int_{B_R(x_0)} |\nabla w|^p  +  C  \int_{B_R(x_0)} |\nabla v|^p \nn\\
&\leq&  C \left (\frac{r}{R} \right)^m  \int_{B_R(x_0)} |\nabla u|^p  +  C  \int_{B_R(x_0)} |\nabla v|^p \nn\\
&\leq&   C  \left (\frac{r}{R} \right)^m    \int_{B_R(x_0)} |\nabla u|^p  +   C \  R^{m-p} (\epsilon_{m,p})^p \ \ ||  \nabla u ||^p_{M^p_p(B_{2R}(x_0)) }.
\eee

For the rest of the proof, we can apply the same arguments as in Schikorra (\cite{Sc}, p.510-511), similarly to Rivi\`{e}re-Struwe (\cite{RS}, Proof of Theorem 1.1, p.459-460),  to obtain the Morrey type estimates for $\nabla u$. The H\"{o}lder continuity of $u$ in $B$ follows immediately from an iteration argument as in  \cite{Gi}.     \eop

\begin{rem} By slightly modifying the proof, we will see that the regularity result
in Theorem 1.5 still hold if the elliptic system \eqref{1.15} is replaced by the following:
\bee    -\ {\rm div}\  \left ( Q \ \nabla u + \Omega  \  u  \right ) = 0, \quad {\rm in}\ \mathcal{D}'(B)
\eee
with $ Q \in W^{1,2}\cap L^{\infty}(B, {\rm GL}(n))$ satisfying $|Q|+|Q^{-1}|\leq \Lambda,  {\rm a.e.\ in}\ B$, for some uniform constant $\Lambda>0$. The proof relies on applying Hodge decomposition to $Q \ \nabla u$ to get the Morrey type estimates for $\nabla u$ as is done by Rivi\`{e}re-Struwe {\rm (}\cite{RS}. See also Schikorra \cite{Sc}{\rm)}.
\end{rem}

\vskip 1.5cm

\section{Harmonic maps into Standard Stationary Lorentzian manifolds}
\vskip0.5cm

In this section, we shall first show that the Euler-Lagrangian equations for weakly harmonic maps into standard stationary Lorentzian manifolds are elliptic systems of the form \eqref{1.3} and then apply Theorem 1.2 to prove the $\epsilon$-regularity (Theorem 1.3) for such maps.

\vskip0.2cm

\noindent{\bf Proof of Theorem 1.3:} Let $(t,u) \in W^{1,2}(B,\mathbb{R}\times M)$ be a weakly harmonic map from
 $B$ into $(\mathbb{R}\times M,g)$, where the metric $g$ is defined as in \eqref{1.9}. For any  $s \in W^{1,2}_0\cap L^{\infty}(B, \mathbb{R}) $ and for any $v \in W^{1,2}_0\cap L^{\infty}(B, \mathbb{R}^n ) $, we have that
\bee t_{\e} = t + \e s \quad {\rm and} \quad u_{\e} = \Pi ( u + \e v)
\eee
are well defined for sufficiently small $\e>0$. Hence $(t_\e, u_\e )\in  W^{1,2}(B,\mathbb{R}\times M) $ gives an admissible variation for $(t,u)$.
By Definition 1.1, there holds
\bee \label{3.2} \left. \frac{d}{d \e} E(t_\e,u_\e) \right|_{\e=0} = 0, \quad \forall s \in W^{1,2}_0\cap L^{\infty}(B, \mathbb{R}), \forall v \in W^{1,2}_0\cap L^{\infty}(B, \mathbb{R}^n )
\eee
A straightforward calculation gives
\bee  \label{3.3} \int_B \left  \{ - \frac{1}{2}(\nabla \b \cdot w) |\nabla t + \omega _i\nabla u^i|^2
-\b (\nabla t + \omega_i\nabla u^i) \cdot \left (\nabla s + \omega_j\nabla w^j+ (\nabla \omega_k \cdot w)\nabla u^k \right ) + \nabla u \cdot \nabla w \right \}=0,
\eee
where $w=d \Pi(u) v, v \in W^{1,2}_0\cap L^{\infty}(B, \mathbb{R}^n )$.

To deduce the Euler-Lagrangian equations, we shall choose appropriate admissible variations in \eqref{3.3}.

First, taking $ s \in W^{1,2}_0\cap L^{\infty}(B)$ and $v \equiv0$ in \eqref{3.3}, we obtain
\bee 0= \int_B  -\b(u)(\nabla t + \omega_i(u)\nabla v^i) \cdot \nabla s. \nn
\eee
Since $s \in W^{1,2}_0\cap L^{\infty}(B)$ is arbitrarily chosen, we get the following conservation law
\bee \label{3.4}   -\ {\rm div} \left \{ \ \beta(u) \ (\nabla t + \omega_i(u)\nabla u^i) \ \right \} = 0.
\eee

Next, taking  $w=d \Pi(u) v, v \in W^{1,2}_0\cap L^{\infty}(B, \mathbb{R}^n)$ and $s \equiv0$ in \eqref{3.3} gives
\bee  \label{3.5}  0 &=& \int_B \left  \{ - \frac{1}{2}(\nabla \b \cdot w) |\nabla t + \omega _i\nabla u^i|^2 -\b (\nabla t + \omega _i\nabla u^i) \cdot \left( \omega_j\nabla w^j+ (\nabla \omega_k \cdot w)\nabla u^k \right ) + \nabla u \cdot \nabla w \right \}\nn\\
&=& \int_B \left  \{ - \frac{1}{2}\left (\frac{\partial \b}{\partial y^j} \cdot w^j \right ) |\nabla t + \omega _i\nabla u^i|^2
-\b (\nabla t + \omega _i\nabla u^i) \cdot\left ( \omega_j\nabla w^j + \nabla u^k \frac{\partial \omega_k}{\partial y^j} \cdot w^j \right )
+ \nabla u \cdot \nabla w \right \} \nn\\
&=& \int_B \left  \{ - {\rm div} \ (\nabla u^j )\cdot w^j + \b(\nabla t + \omega_i\nabla u^i) \cdot  \nabla u^k  \left (   \frac{\partial \omega_j}{\partial y^k} -  \frac{\partial \omega_k}{\partial y^j} \right)  \cdot w^j  - \frac{1}{2}|\nabla t + \omega _i\nabla u^i|^2 \frac{\partial \b}{\partial y^j} \cdot w^j  \right \}
\eee
where in the last step we have used \eqref{3.4} and integration by part. Denote $H:=(H^1,...H^n)$ with
\bee \label{3.6}  H^j : = \beta \left (\nabla t + \omega_i\nabla u^i \right ) \cdot \nabla u^k \left(\frac{\partial \omega _j}{\partial y^k}- \frac{\partial \omega_k}{\partial y^j} \right) - \frac{1}{2}\frac{\partial \b}{\partial y^j}\left |\nabla t + \omega_i\nabla u^i\right |^2.
\eee
Then \eqref{3.5} becomes
\bee 0 &=& \int_B \left ( - \ {\rm div} \ \nabla u  + H  \right ) \cdot d \Pi(u) v, \quad \forall  v \in W^{1,2}_0\cap L^{\infty}(B, \mathbb{R}^n)  \nn
\eee
Since $ v \in W^{1,2}_0\cap L^{\infty}(B, \mathbb{R}^n)$ is arbitrarily chosen, we have (similarly to the calculations in \cite{H2}, Chapter 1.)
\bee  \label{3.7}  - \ {\rm div} \ \nabla u  - A(u)(\nabla u, \nabla u) +  d \Pi(u) H=0 ,
\eee
where $A $ is the second fundamental form of $M $ in $\mathbb{R}^n$. Let $\nu_l, l=d+1,d+2,...,n$ be an orthonormal frame for the normal bundle $T^\bot M$
(and still denote by $\nu_l$ the corresponding normal frame along the map $u$), then we can rewrite \eqref{3.7} as follows:
\bee \label{3.8}  - \ {\rm div} \ \nabla u = \nu_l \nabla \nu_l \cdot \nabla u -  H + \langle H , \nu_l \rangle \nu_l,
\eee
where $\la \cdot,\cdot \ra$ denotes the Euclidean metric on $\mathbb{R}^n$. We have thus obtained the Euler-Lagrangian equations:
\bee \label{3.9}   -\ {\rm div} \left \{\ \beta(u)  \ (\nabla t + \omega_i(u)\nabla u^i) \ \right \} = 0,
\eee
\bee \label{3.10}  -\ {\rm div} \ \nabla u = \nu_l \nabla \nu_l \cdot \nabla u -  H + \langle H , \nu_l \rangle \nu_l.
\eee

To proceed, we write the system of equations \eqref{3.9} and \eqref{3.10} in the form of \eqref{1.3}.
By Hodge decomposition, we conclude from the conservation law \eqref{3.4} that there exists $\eta \in W^{1,2}(B, \wedge^2 \mathbb{R}^{m})$ such that
\bee \label{3.11}  \beta(u) \ (\nabla t + \omega_i(u)\nabla u^i) = {\rm curl}\  \eta.
\eee
Then, by \eqref{3.6}, we can rewrite the equation \eqref{3.10} as:
\bee \label{3.12} -\ {\rm div} \ \nabla u^j  = \Theta^{jk} \cdot \nabla u^k  + a_{jk} \ {\rm curl}\ \eta  \cdot \nabla u^k + b_j \ {\rm curl}\ \eta \cdot \b(u) (\nabla t + \omega_i\nabla u^i),
\eee
where
\bee  \label{3.13} \Theta^{jk}&: =&  \nu_l^j \nabla \nu_l^k -  \nu_l^k \nabla \nu_l^j  \\
\label{3.14} a_{jk}&:=& - \left(\frac{\partial \omega _j}{\partial y^k}- \frac{\partial \omega_k}{\partial y^j} \right)+ \left(\frac{\partial \omega _i}{\partial y^k}
- \frac{\partial \omega_k}{\partial y^i} \right)  \nu_l^i \nu_l^j    \\
\label{3.15} b_j&:=& \frac{1}{2\b^2(u)}\left ( \frac{\partial \b}{\partial y^j}  - \frac{\partial \b}{\partial y^i }  \nu_l^i \nu_l^j  \right )
\eee
Now we can write the Euler-Lagrangian equations \eqref{3.9} and \eqref{3.10} as the following elliptic system:
\bee \label{3.16} -\ {\rm div} \left \{Q \cdot \left ( \ba{ccc}
            \nabla t    \\
           \nabla u   \ea \right)   \right \}  =  \Theta \cdot  Q  \ \left ( \ba{ccc}
            \nabla t    \\
           \nabla u   \ea \right )  + F  \ {\rm curl}\ \zeta \cdot Q \ \left ( \ba{ccc}
            \nabla t    \\
           \nabla u   \ea \right ),
\eee
where
\bee  \label{3.17} Q &=&  \tilde{Q} \circ  u,  \quad
   \tilde{Q}  =   \left ( \ba{ccc}
            \b  & \b\omega    \\
            0   &  I_n   \ea \right ), \quad \omega= (\omega_1,\omega_2,...,\omega_n) \\
\label{3.18} \Theta &=&  \left ( \ba{ccc}
            0  &    0    \\
            0   &  (\Theta^{jk})   \ea \right )  \\
\label{3.19}  F & =&  \left ( \ba{ccc}
           0  & 0 \cdots 0    \\
            b_1     &  a_{11}   \cdots  a_{1n}   \\
                   \vdots               &       \ddots                                                 \\
            b_n    &  a_{n1}   \cdots  a_{nn}  \\   \ea \right )      \\
\label{3.20} \zeta &=&   {\rm diag}\ ( \eta, \eta, ...,\eta).
\eee

Since $M$ is compact, $\b \in C^2(M, (0,\infty))$ and $\omega \in C^2(\Omega^1(M))$,
there exists $\lambda >0$ depending only on the target $(\mathbb{R}\times M,g)$ such that for any $y\in M$ there hold
\bee \label{3.21} 0 < \lambda^{-1} \leq \b(y), \quad |\b(y)|+ | \nabla \b(y)| + | \nabla^2 \b(y)|\leq \lambda < \infty, \quad |\omega (y)| + |\nabla \omega (y)| + |\nabla^2\omega (y)|\leq \lambda< \infty,
\eee
Using the notations \eqref{3.13}-\eqref{3.15}, \eqref{3.17}-\eqref{3.20} and the above estimates \eqref{3.21}, we can easily verify that $\Theta \in L^{2}(B,so(n+1)\otimes \wedge ^1 \mathbb{R}^{m})$,
 $F \in W^{1,2}\cap L^{\infty}(B, {\rm M}(n+1) )$, $\zeta \in W^{1,2}(B, {\rm M}(n+1) \otimes \wedge^2 \mathbb{R}^{m})$,
  $ Q \in W^{1,2}\cap L^{\infty}(B, {\rm GL}(n+1))$ and the following estimates hold:
\bee \label{3.22} |Q|+|F| \leq C_1(\lambda), \quad {\rm a.e.\ in}\ B
\eee
and
\bee  \label{3.23}     || \Theta  ||_{M^2_2(B)} +   || \nabla F ||_{M^2_2(B)} +  || \nabla  Q  ||_{M^2_2(B)}  \leq  C_2(\lambda) \   || \nabla u  ||_{M^2_2(B)},
\eee
where $C_1(\lambda)>0$ and $C_2(\lambda)>0$ are constants also depending on $\lambda$.

To estimate $|Q^{-1}|$, we note that
\bee  \tilde{Q}^{-1}  =   \left ( \ba{ccc}
            \b^{-1}  & -\omega    \\
            0   &  I_n   \ea \right ). \nn
\eee
Hence, by \eqref{3.21}, there exists some constant $C_3(\lambda) > 0$ such that
\bee \label{3.24} |Q^{-1}| =  |\tilde{Q}^{-1} \circ  u|  \leq C_3(\lambda), \quad {\rm a.e.\ in}\ B.
\eee

On the other hand, it follows from \eqref{3.11} and \eqref{3.21} that
\bee | {\rm curl}\  \eta |  \leq C_4(\lambda)  \left ( | \nabla t | +  | \nabla u| \right ), \quad {\rm a.e.\ in} \ B.\nn
\eee
By \eqref{3.20} and the above inequality, we verify that
\bee \label{3.25}  ||{\rm curl}\  \zeta ||_{M^2_2 (B)} \leq C_5(\lambda)  \  \left (   || \nabla t  ||_{M^2_2(B)} +  || \nabla u  ||_{M^2_2(B)} \right).
\eee

Combining \eqref{3.22} and \eqref{3.24} gives
\bee \label{3.26} |Q|+|Q^{-1}|+|F| \leq C_1(\lambda)+C_3(\lambda), \quad {\rm a.e.\ in}\ B.
\eee
Combining \eqref{3.23} and \eqref{3.25} gives
\bee \label{3.27} && \hskip0.3cm || \nabla t ||_{M^2_2(B)} +  || \nabla u  ||_{M^2_2(B)} +  ||{\rm curl}\  \zeta ||_{M^2_2 (B)} +  || \Theta  ||_{M^2_2(B)} +   || \nabla F ||_{M^2_2(B)} +  || \nabla  Q  ||_{M^2_2(B)} \nn \\
&& \leq (1+  C_2(\lambda) + C_5(\lambda))  \  \left (|| \nabla t  ||_{M^2_2(B)} +  || \nabla u  ||_{M^2_2(B)} \right).
\eee

Take $\Lambda:= C_1(\lambda) + C_3(\lambda)>0$, then $\Lambda$ depends only on $(\mathbb{R}\times M,g)$.
Let $\epsilon_{m,\Lambda}>0$ be the small constant (depending on $m$ and $\Lambda$) as in Theorem 1.2. Take $$\epsilon_{m}:=\frac{\epsilon_{m,\Lambda}}{1+  C_2(\lambda) + C_5(\lambda)},$$ then $\epsilon_{m}>0$ depends only on $(\mathbb{R}\times M,g)$. Applying Theorem 1.2 to the elliptic system \eqref{3.16}, we conclude from \eqref{3.26} and \eqref{3.27} that $(t,u)$ is H\"{o}lder continuous in $B$ if $|| \nabla t  ||_{M^2_2(B)} +  || \nabla u  ||_{M^2_2(B)} < \epsilon_{m}.$ By standard elliptic regularity theory, $(t,u)$ is as smooth as the regularity of the target $(\mathbb{R}\times M,g)$ permits.   \hfill $\Box$

\vskip1.5cm

\section{Harmonic maps into pseudospheres $\mathbb{S}^n_\nu$ $(1\leq \nu \leq n)$}
\vskip0.5cm

In this section, we shall first prove Proposition 1.1 and Proposition 1.2. Then, with the help of these two propositions, we apply Theorem 1.5 to prove the regularity results (Theorem 1.6 and Theorem 1.7) for weakly harmonic maps into pseudospheres $\mathbb{S}^n_\nu  \ (1 \leq \nu \leq n)$.

\vskip0.2cm

\noi{\bf Proof of Proposition 1.1:} Fix  $i\neq j \in \{1,2,...,n+1\}$. Let $E_{ij}\in so(n+1)$ be the matrix
whose $(i,j)$-component is $1$, $(j,i)$-component is $-1$ and all the other components are $0$. Let $\mathcal{E}$
 be the matrix defined as in \eqref{1.17}. Then one verifies that $E_{ij}\mathcal{E} \in so( \nu, n+1-\nu )$ and
$e^{E_{ij}\mathcal {E}} \in  {\rm O}(\nu, n+1-\nu )$ (see e.g. \cite{O}). For any $\varphi \in C_0^{\infty}(B)$, define
\bee R_t:=e^{t\varphi E_{ij}\mathcal {E}} \in  C_0^{\infty}(B, {\rm O}(\nu, n+1-\nu )).
\eee
Using the property of an element in the group ${\rm O}(\nu, n+1-\nu )$ (see \eqref{1.19}), we have
\bee \label{4.2} \la R_t u, R_t u \ra_{\mathbb{R}^{n+1}_\nu} = (R_t u)^T \mathcal {E} \ R_t u = u^T R_t ^T \mathcal {E} \ R_t u = u^T  \mathcal {E} \ u =1, \quad {\rm a.e.\  in}\  B.
\eee
It follows that $R_t u \in W^{1,2}(B,\mathbb{S}^n_\nu )$. Since $u$ is weakly harmonic, by Definition 1.2, we calculate
\bee \label{4.3} 0= \left. \frac{d}{dt} \right|_{t=0} E(R_t u)  &=& \int_B    (\nabla (R_{0} u))^T \mathcal {E}\  \left. \frac{d}{dt} \right|_{t=0}\ (\nabla (R_t u)) \nn\\
&=& \int_B (\nabla  u)^T \mathcal {E}\ \left(\nabla ( \varphi E_{ij}\ \mathcal {E}\ u ) \right)      \nn\\
&=& \int_B (\nabla  u)^T \mathcal {E}\ E_{ij}\ \mathcal {E}\ u \  \nabla \varphi + (\nabla  u)^T \mathcal {E} \ E_{ij}\ \mathcal {E}\ \nabla u \  \varphi \nn\\
&=& \int_B (\nabla  u)^T \mathcal {E}\ E_{ij}\ \mathcal {E}\ u  \  \nabla \varphi \nn\\
&=& (\varepsilon_{ii} \varepsilon_{jj})  \int_B  \left(u^i \nabla u^j - u^j \nabla u^i \right ) \nabla \varphi,
\eee
where we have used the fact that $\mathcal {E}\ E_{ij}\ \mathcal {E} \in so(n+1)$
and hence
\bee (\nabla  u)^T \mathcal {E}\ E_{ij}\ \mathcal {E}\ \nabla u = 0 \quad {\rm a.e.\  in}\  B.\nn
\eee
Since $\varphi \in C_0^{\infty}(B)$ is arbitrary and $\varepsilon_{ii} \varepsilon_{jj}$ is either $1$ or $-1$ (see \eqref{1.17}), we conclude from \eqref{4.3} that the conservation laws \eqref{1.27} hold for $i\neq j$.

The case of $i=j$ is trivial. This completes the proof.  \eop

\vskip 0.2cm

\noi{\bf Proof of Proposition 1.2:} First, by definition of the space
$W^{1,2}(B,\mathbb{S}^n_\nu )$ (see \eqref{1.21}), we have
\bee \label{4.4}  u^j  \varepsilon_{jk}  u^k = 1 \quad {\rm a.e.\  in}\  B.
\eee
Taking $\nabla$ on both sides of \eqref{4.4} gives
\bee \label{4.5} \nabla  u^j  \varepsilon_{jk}  u^k = 0 \quad {\rm a.e.\  in}\  B.
\eee
Recall that (see \eqref{1.28}) $\Theta = \left (\Theta^{ij} \right ) = \left (  u^i \nabla u^j - u^j \nabla u^i \right )$. Combining \eqref{4.4} and \eqref{4.5}, we calculate
\bee \label{4.6} \nabla u^i + \Theta^{ij}  \varepsilon_{jk}  u^k & = &  \nabla u^i + \left (  u^i \nabla u^j - u^j \nabla u^i \right )  \varepsilon_{jk}  u^k  \nn\\
&=&  \nabla u^i \left (1-   u^j  \varepsilon_{jk}  u^k  \right ) + u^i \left (  \nabla  u^j  \varepsilon_{jk}  u^k\right )   \nn\\
&=& 0 \hskip 5cm {\rm a.e.\  in}\  B.
\eee
This proves \eqref{1.29}.

Since $u\in W^{1,2}(B, \mathbb{R}^{n+1})$, one verifies that $\nabla u^i + \Theta^{ij} \varepsilon_{jk} u^k \in L^1(B)$ for each $i$. Taking $ - \ {\rm div}$ on both sides of \eqref{4.6} gives
\bee   - \ {\rm div} \left ( \nabla u + \Theta \ \mathcal {E}\  u \right ) = 0,  \quad {\rm in}\ \mathcal{D}'(B). \nn
\eee

Next, we assume that $u$ is weakly harmonic and for any fixed $1<p< \frac{m}{m-1}$ there holds $||\nabla u ||_{ M^p_p (B)} <  \infty$. We shall derive the estimate \eqref{1.31}.

Let $q = \frac{p}{p-1}>m$ be the conjugate exponent of $p$. Let $B_R(x_0) \subset B_{1/2}$. For any $\Phi\in L^q(B_R(x_0), \wedge ^1 \mathbb{R}^{m} )$ with $||\Phi||_{ L^q(B_R(x_0))} \leq 1$ and for any $0<\rho< R$, let $\tau = \tau(\rho)\in C_0^{\infty}(B_R(x_0), [0,1]) $ be a cut-off function satisfying $$\tau \equiv 1, \quad {\rm on}\ B_\rho(x_0),$$ then $ \tau  \Phi$ is supported in $B_R(x_0)$ and vanishes on $\partial B_R(x_0)$. By Hodge decomposition, there exist $\a \in W^{1,q}_0(B_R(x_0))$, $\b \in W^{1,q}_0(B_R(x_0), \wedge ^2 \mathbb{R}^{m})$ and a harmonic $h\in C^{\infty}(B_R(x_0), \wedge ^1 \mathbb{R}^{m} )$ such that
\bee   \label{4.7}   \tau   \Phi = \nabla \a + {\rm curl}\ \b+h.
\eee
Moreover, we have
\bee  \label{4.8} ||  \nabla \a || _{ L^q(B_R(x_0))}   +    ||  \nabla  \b || _{ L^q(B_R(x_0))}  \leq  C  ||  \tau \Phi || _{ L^q(B_R(x_0))}    \leq  C   ||  \Phi || _{ L^q(B_R(x_0))} \leq C,
\eee
where $C>0$ is a constant independent of $\rho$ and $R$. Recall that $ \tau \in C_0^{\infty}(B_R(x_0))$, we get $ \left. h \ \right|_{\partial B_R(x_0)} =\left. (\tau  \Phi) \right|_{\partial B_R(x_0)}=0$. Since $h$ is harmonic, it follows that $h\equiv 0$ in $B_R(x_0)$.

Since $u$ is weakly harmonic, by Proposition 1.1, $\Theta = \left (\Theta^{ij} \right ) = \left (  u^i \nabla u^j - u^j \nabla u^i \right )$ is divergence free. Then, using \eqref{4.7}, \eqref{4.8} and the fact that $h\equiv 0$ in $B_R(x_0)$, and applying Lemma 2.1, we estimate for fixed $i,j \in \{1,2,...,n+1\}$,
\bee \int_{B_R(x_0)} \left (\tau \Theta^{ij} \right ) \cdot \Phi & =&  \int_{B_R(x_0)}  \Theta^{ij} \cdot  (\tau \Phi) \nn\\
&=&  \int_{B_R(x_0)}  \Theta^{ij} \cdot  \left (\nabla \a + {\rm curl}\ \b \right ) \nn\\
&=&  \int_{B_R(x_0)}  \Theta^{ij}  \cdot  {\rm curl}\ \b  \nn\\
&=&  \int_{B_R(x_0)}  \left ( u^i \nabla u^j - u^j \nabla u^i \right ) \cdot {\rm curl}\ \b \nn\\
&=&  \int_{B_R(x_0)} \left \{ \left (\nabla u^j  \cdot {\rm curl}\ \b \right ) u^i   -  \left (\nabla u^i \cdot {\rm curl}\ \b \right ) u^j \right\} \nn\\
&\leq& C \ ||   \nabla u || _{ L^p(B_R(x_0))} \ ||  {\rm curl}\ \b|| _{ L^q(B_R(x_0))} \ ||  \nabla u || _{M^p_p(B_{2R}(x_0)) }   \nn\\
&\leq& C \ ||   \nabla u || _{ L^p(B_R(x_0))} \ ||  \nabla  \b || _{ L^q(B_R(x_0))} \ ||  \nabla u || _{M^p_p(B_{2R}(x_0)) }   \nn\\
&\leq&  C\  ||  \nabla u || _{ L^p(B_R(x_0))}  \ ||  \nabla u || _{M^p_p(B_{2R}(x_0)) }
\eee
By duality characterization of $L^p$ functions, we have
\bee  ||  (\tau \Theta^{ij}) || _{ L^p(B_R(x_0))}  \leq  C\   ||   \nabla u || _{ L^p(B_R(x_0))}  \ ||  \nabla u || _{M^p_p(B_{2R}(x_0)) }.
\eee
It follows that
\bee  || \Theta^{ij}|| _{ L^p(B_\rho(x_0))} \leq  ||  (\tau \Theta^{ij}) || _{ L^p(B_R(x_0))}  \leq C\  ||  \nabla u || _{ L^p(B_R(x_0))}  \ ||  \nabla u || _{M^p_p(B_{2R}(x_0)) }.
\eee
Since $\rho \in (0, R)$ is arbitrary, let $\rho \nearrow R$, then we get
\bee  || \Theta^{ij}|| _{ L^p(B_R(x_0))}  \leq  C\   ||  \nabla u || _{ L^p(B_R(x_0))}  \ ||  \nabla u || _{M^p_p(B_{2R}(x_0)) }.
\eee
Furthermore, using the definition of the Morrey norm  $||  \nabla u || _{ M^p_p(B)}$ and the fact that $B_{2R}(x_0) \subset B$, we estimate
\bee || \Theta||_{ L^p(B_R(x_0))} = \sum_{i,j}|| \Theta^{ij}||_{ L^p(B_R(x_0))} &\leq& C\   ||\nabla u ||_{ L^p(B_R(x_0))} \ ||  \nabla u || _{M^p_p(B_{2R}(x_0)) }  \nn\\
&\leq&  C\  R^{\frac{m}{p}-1} \  ||  \nabla u || _{ M^p_p(B)}  \  ||  \nabla u || _{ M^p_p(B)} \nn\\
&=&  C\   R^{\frac{m}{p}-1} \  ||  \nabla u ||^2 _{ M^p_p(B)}. \nn
\eee
Since the ball $B_R(x_0) \subset B_{1/2}$ is arbitrary, it follows that
\bee
|| \Theta ||_{M^p_p(B_{1/2})} = \underset {B_R(x_0)\subset B_{1/2}}{\rm sup}  \left (  R^{p-m} \int_{B_R(x_0)} |\Theta|^p \right )^{\frac{1}{p}} \leq  C \  ||  \nabla u ||^2_{ M^p_p(B)}. \nn
\eee
Thus, we have completed the proof. \eop

\vskip0.2cm

\noindent{\bf Proof of Theorem 1.6:} Note that $\mathcal{E}$ is a constant matrix. Combining Proposition 1.1, Proposition 1.2, Theorem 1.5 and using a rescaling of the domain gives that $u$ is  H\"{o}lder continuous in $B$. Moreover, since ${\rm div}\ \Theta  = 0 $,
we can rewrite the equation in \eqref{1.30} as
\bee \label{}  -\ {\rm div}\ \nabla u =  \Theta \ \mathcal {E} \cdot \nabla  u. \nn
\eee
By standard elliptic regularity theory, $u$ is smooth in $B$. \eop

\vskip0.2cm

\noindent{\bf Proof of Theorem 1.7:} Fix some  $1<p< \frac{m}{m-1}=2$. By conformal invariance in dimension $m=2$ and rescaling in the domain, we assume W.L.O.G that
\bee
  || \nabla u ||^2_{ L^2 (B)}  < \epsilon_{2, p},
\eee
where  $ \epsilon_{2,p}$ is given in Theorem 1.6 with $m=2$. By a straightforward calculation, it follows that
\bee
 || \nabla u ||^2_{ M^p_p (B)} \leq  || \nabla u ||^2_{ L^2 (B)}   <  \epsilon_{2, p}.
\eee
Applying Theorem 1.6  with $m=2$ gives that $u$ is H\"{o}lder continuous (and hence smooth) in $B$. \eop

\vskip1.5cm

\section{Generalized (weakly) harmonic maps into $\mathbb{S}^n_\nu$ $(1\leq \nu \leq n)$}
\vskip0.5cm

In this section, we shall prove the $\epsilon$-regularity result (Theorem 1.8) for generalized (weakly) harmonic maps into $\mathbb{S}^n_\nu$ ($1\leq \nu \leq n$). Throughout this section,  $B$ will denote the unit disc in $\mathbb{R}^2 $.

\vskip0.2cm

\noindent{\bf Proof of Theorem 1.8:} Slightly modifying some arguments in the proofs of Proposition 1.2
and Theorem 1.5 will be sufficient to prove this theorem.

Fix any $\frac{4}{3}<p<2$ and let $u\in W^{1,\frac{4}{3}}(B,\mathbb{S}^n_\nu)$ $(1\leq \nu \leq n)$ be a generalized
(weakly) harmonic map  satisfying
\bee \label{5.1}
 || \nabla u ||^2_{ M^p_p (B)} < \epsilon_{p}
\eee
with $\epsilon_{p}>0$ being determined later. Then $u \in W^{1,p}(B)$ and hence $\Theta = ( \Theta^{ij} ) : = \left( u^i \nabla u^j - u^j \nabla u^i \right) \in L^{p'}(B)$,
 where $p'=\frac{2p}{4-p} \in (1,p)$. By Definition 1.5, there holds
\bee  \label{5.2}  {\rm div } \ \Theta= 0, \quad {\rm in}\ \mathcal{D}'(B).
\eee
Applying similar arguments as in the proof of Proposition 1.2 (with $m=2$) gives that
\bee \label{5.3} \nabla u + \Theta \  \mathcal{E} \  u =  0 \hskip 1cm {\rm a.e.\  in}\  B
\eee
and
\bee \label{5.4}
|| \Theta ||_{M^{p'}_{p'}(B_{1/2})} \leq C_{p'}\ || \nabla u ||^2_{ M^{p'}_{p'} (B)  }  \leq  C_{p'}\  || \nabla u ||^2_{ M^{p}_{p} (B)  }  \leq  C_{p'}\      \epsilon_{p}.
\eee

Let $B_{2R}(x_0) \subset B_{1/2}$ and let $w \in W^{1,p'}(B_{R}(x_0),\mathbb{R}^{n+1}) $ be solving
\bee \label{5.5} \left \{    \ba{rcll}
  - \ {\rm div}\ \nabla w &=& 0 ,  &\qquad\text{in } B_R(x_0)    \\
               w &=& u,        &\qquad\text{on }   \partial B_R(x_0) \ea \right.
\eee
and define $v:=u-w \in W^{1,p'}_{0}(B_R(x_0),\mathbb{R}^{n+1})$.

Let $q' = \frac{p'}{p'-1}$ be the conjugate exponent of $p'$. Then for any $\varphi\in W^{1,q'}_0(B_R(x_0))$ with $||\varphi||_{ W^{1,q'}(B_R(x_0))} \leq 1$.
Using \eqref{5.3} and \eqref{5.5}, we get
\bee   \int_{B_R(x_0)}  \nabla v^i \cdot \nabla \varphi =  \int_{B_R(x_0)} \nabla u^i \cdot \nabla \varphi - \int_{B_R(x_0)} \nabla w^i \cdot \nabla \varphi =  - \  \varepsilon_{jj} \int_{B_R(x_0)}  \Theta^{ij}\ u ^j \cdot  \nabla \varphi  \nn
\eee
Then using \eqref{5.1}, \eqref{5.2}, \eqref{5.4}, Lemma 2.1 and taking $\epsilon_{p}>0$ sufficiently small, we can apply the same arguments as in the proof of Theorem 1.5 (with $m=2$) and use a rescaling of the domain to conclude that $u$ is H\"{o}lder continuous and hence smooth (by standard elliptic regularity) in $B$. \eop

\vskip0.2cm

Furthermore, we observe that the $\epsilon$-regularity result in Theorem 1.8 still hold if the Morrey norm  $||\nabla u ||_{ M^p_p (B)}$  is replaced with the Lorentz norm  $|| \nabla u ||_{ L^{(2,\infty)}(B)}$ (which was used in Almeida \cite{Al}). To see this, we recall the following:

\begin{lem}[Almeida \cite{Al}, Lemma 9]  Suppose $D$ has finite measure. Let $1<p<p_1<\infty$.
 Then, there is a constant $C$ such that, for all $q,q_1 \in [1, \infty]$ and for any $f\in L^{(p_1,q_1)}(D)$,
\bee ||f||_{L^{(p,q)}} \leq C \ (\mu(D))^{\frac{p_1-p}{pp_1}} \ ||f||_{L^{(p_1,q_1)}}
\eee
\end{lem}

Recall that $L^{(p,p)}=L^p$. Consequently, we have

\begin{lem}   Let $1<p<2$. Then, there is a constant $C$ such that, for any $f\in L^{(2,\infty)}(B)$,
\bee ||f||_{M^p_p(B)} \leq C \ ||f||_{L^{(2,\infty)}(B)}
\eee
\end{lem}

\pr  Take $p=q$, $p_1=2$, $q_1=\infty$ in Lemma 5.1 and let $D$ run over all discs $B_R(x_0)\subset B$. \eop

\vskip0.2cm

Combining Theorem 1.8 and Lemma 5.2 gives the following $\epsilon$-regularity result (using the Lorentz norm).

\begin{thm}  There exists $\epsilon>0 $ such that any generalized (weakly) harmonic map $u\in W^{1,\frac{4}{3}}(B, \mathbb{S}^n_\nu)$ $(1\leq \nu \leq n )$ satisfying
\bee
 || \nabla u ||_{ L^{(2,\infty)}(B)} < \epsilon
\eee
is smooth  in $B$.
\end{thm}

\vskip1.5cm

\section{Regularity for an elliptic system with a potential in $so(1,1)$ }
\vskip0.5cm

Throughout this section,  $B$ will denote the unit disc in $\mathbb{R}^2 $. We consider the elliptic system \eqref{1.1} with
 a potential $\Omega\in  L^{2}(B,so(1,1)\otimes \wedge ^1 \mathbb{R}^{2})$. By Hodge decomposition, there exist $\Omega_1 \in  W^{1,2}(B,so(1,1) )$ and $\Omega_2 \in W^{1,2}(B,so(1,1)\otimes \wedge^2 \mathbb{R}^{2} ) $ such that
\bee \label{6.1} \Omega = \nabla \Omega_1 +  {\rm curl} \ \Omega_2,
\eee

\begin{thm} Let $u\in W^{1,2}(B,\R^2)$ be a weak solution of the elliptic system \eqref{1.1} with  a potential $\Omega\in  L^{2}(B,so(1,1)\otimes \wedge ^1 \mathbb{R}^{2})$. Decompose $\Omega$ as in \eqref{6.1}. If $\Omega_1 \in L^{\infty}(B, so(1,1))$, then $u$ is  H\"{o}lder continuous in $B$.
\end{thm}

\pr Since $\Omega_1$ takes values in $so(1,1)$, we can write (see O'Neill's book \cite{O})
\bee \label{6.2}  \Omega_1 = \left ( \ba{ccc}
            0 & s   \\
            s &  0  \ea \right ), \quad {\rm for \ some}\ s  \in W^{1,2}(B).
\eee
Consequently, we have $\nabla \Omega_1  \Omega_1 = \Omega_1  \nabla \Omega_1$ and hence  $\nabla ( e^{\Omega_1} ) = e^{\Omega_1} \nabla \Omega_1$. Then we calculate
\bee \label{6.3} - \ {\rm div} \left ( e^{\Omega_1} \nabla u \right ) = - \ e^{\Omega_1}  \nabla \Omega_1 \cdot \nabla u + e^{\Omega_1}  \Omega \cdot \nabla u = e^{\Omega_1} {\rm curl} \ \Omega_2 \cdot \nabla u.
\eee
Using \eqref{6.2}, we get
\bee e^{\Omega_1} = (e^{\Omega_1})^T =  \left ( \ba{ccc}
             {\rm cosh}\ s &   {\rm sinh}\  s   \\
            {\rm sinh}\  s &   {\rm cosh}\  s  \ea \right ), \quad
 e^{-\Omega_1}=(e^{\Omega_1})^{-1} = \left ( \ba{ccc}
             {\rm cosh}\  s   &   -  {\rm sinh}\  s   \\
            -  {\rm sinh}\ s &   {\rm cosh}\  s  \ea \right ). \nn
\eee

Since $\Omega_1 \in L^{\infty}(B)$, there exists a constant $\lambda \in (0,\infty)$, such that $ |s| \leq \lambda, {\rm a.e. \ in }\  B$.
Therefore, we have
$$|e^{\Omega_1}| +  |(e^{\Omega_1})^{-1}|  \leq C(\lambda), \quad  {\rm a.e. \ in }\  B $$
for some constant $C(\lambda)>0$ depending on $\lambda$.

On the other hand, one verifies that $e^{\Omega_1} \in W^{1,2}\cap L^{\infty}(B, {\rm M}(2))$. Recall that $\Omega_2 \in W^{1,2}(B, so(1,1)\otimes \wedge^2 \mathbb{R}^{2})$. Applying Theorem 1.2 (with $m=2$ and $\Lambda=C(\lambda)$) to the elliptic system \eqref{6.3}, using the conformal invariance in dimension $m=2$ and rescaling in the domain, we get the H\"{o}lder continuity of $u$ in $B$.    \eop

\vskip0.2cm

Theorem 6.1 is optimal. To see this, we set
\bee  s (x) = {\rm  log\  log} \ \frac{2} {|x|}, \quad u_1 (x) = {\rm log\  log} \ \frac{2} {|x|},  \quad u_2 (x) = {\rm  log\  log} \ \frac{2} {|x|}, \quad x \in B.  \nn
\eee
Then the map $u=(u_1,u_2)^T \in W^{1,2}(B, \mathbb{R}^{2}) $ is a weak solution to the elliptic system \eqref{1.1} with a potential $\Omega$ satisfying
\bee \Omega = \left ( \ba{ccc}
            0 & \nabla s   \\
             \nabla s  &  0  \ea \right ) \in   L^{2}(B,so(1,1)\otimes \wedge ^1 \mathbb{R}^{2}) \ {\rm and }\  s \ {\rm is\ not\ in} \  L^{\infty}(B). \nn
\eee
However, $u$ is not in $L^{\infty}(B)$.

\end{document}